\documentclass{article}
\usepackage{amssymb,amsmath}
\usepackage{graphics}

\def\qed{\hfill {\hbox{${\vcenter{\vbox{               
   \hrule height 0.4pt\hbox{\vrule width 0.4pt height 6pt
   \kern5pt\vrule width 0.4pt}\hrule height 0.4pt}}}$}}}

\def\tr{\triangleright}
\def\bar{\overline}

\newtheorem{definition}{Definition}

\newtheorem{example}{Example}

{\qed\par\medskip}

\date{}

\title{\Large \textbf{Generalized quandle polynomials}}

\author{Sam Nelson}

\begin{document}
\maketitle

\begin{abstract}
We define a family of generalizations of the two-variable quandle polynomial.
These polynomial invariants generalize in a natural way to eight-variable
polynomial invariants of finite biquandles. We use these polynomials to define
a family of link invariants which further generalize the quandle counting 
invariant.
\end{abstract}



\section{\large \textbf{Introduction}}

In \cite{N} a two-variable polynomial invariant of finite quandles, denoted
$qp_Q(s,t)$,
was introduced. This invariant was shown to distinguish all non-Latin 
quandles of order 5 and lower. A slight modification gives an invariant
of subquandles embedded in larger quandles which is capable of 
distinguishing isomorphic subquandles embedded in different ways. 
This subquandle polynomial was used to augment the quandle
counting invariant $|\mathrm{Hom}(Q(L),T)|$ to obtain a multiset-valued
invariant which can distinguish knots and links with the same quandle
counting invariant value.

In this paper we generalize the quandle polynomial in two ways. In section 
\ref{gqp} we define a family of two-variable 
polynomial invariants of finite quandles indexed by pairs of integers, 
denoted $qp_{m,n}(Q)$, such that $qp_Q(s,t)=qp_{1,1}(Q)$. In section 
\ref{bqp} we extend these generalized quandle polynomials in a natural 
way to obtain a family of eight-variable polynomial invariants of finite 
biquandles indexed by pairs of integers, denoted $bp_{m,n}(B)$. In section 
\ref{inv} we define and give examples of link invariants defined using 
generalized quandle polynomials. 
In section \ref{questions} we collect a few questions for future research.

\section{\large \textbf{The $(m,n)$ quandle polynomial}}\label{gqp}

We begin with a definition from \cite{J}.

\begin{definition}
\textup{A \textit{quandle} is a set $Q$ with a binary operation 
$\tr:Q\times Q\to Q$ satisfying}
\begin{list}{}{}
\item[(i)]{\textup{for every $a\in Q$, $a\tr a=a$,}}
\item[(ii)]{\textup{for every $a,b\in Q$, there is a unique $c\in Q$ 
such that $a=c\tr b$, and}}
\item[(iii)]{\textup{for every $a,b,c\in Q$ we have 
$(a\tr b)\tr c=(a\tr c)\tr (b\tr c)$.}}
\end{list}
\textup{Axiom $(ii)$ is equivalent to}
\begin{list}{}{}
\item[(ii')]{\textup{There is a second operation $\tr^{-1}:Q\times Q\to Q$
such that $(a\tr b)\tr^{-1}b=a=(a\tr^{-1}b)\tr b$.}}
\end{list}
\end{definition}

As a useful abbreviation, let us denote the $n$-times repeated quandle
operation as 
\[x\tr^n y= (\dots (x\tr y) \tr y \dots ) \tr y\]
and
\[x\tr^{-n} y= (\dots (x\tr^{-1} y) \tr^{-1} y\dots )\tr^{-1} y\]
where as expected, $n$ is the number of triangles. Note that $x\tr^0 y=x$ 
for all $x,y$.

A standard example of a quandle is any group $G$, which has quandle
structures given by 
$g\tr h= h^{-n}gh^n$ for $n\in \mathbb{Z}$ and $g\tr h= t(gh^{-1})h$
for any $t\in\mathrm{Aut}(G),\ g,h\in G$. The special case of the latter
where $G$ is abelian is called an \textit{Alexander quandle} and may
be regarded as a module over $\mathbb{Z}[t^{\pm 1}]$ by thinking of 
$t\in \mathrm{Aut}(G)$ as a formal variable; in additive notation 
we have \[x\tr y = tx + (1-t)y.\]

Another standard example of a quandle structure is any module $V$ over 
a commutative ring $R$ with an antisymmetric bilinear form\footnote{If the 
characteristic of $R$ is 2, 
then we require that $\langle\mathbf{x},\mathbf{x}\rangle=0$.} 
$\langle,\rangle:V\times V\to R$ with $\mathbf{x}\tr \mathbf{y} 
= \mathbf{x}+\langle\mathbf{x},\mathbf{y}\rangle\mathbf{y}$.
When $R$ is a field and $\langle,\rangle$ is nondegenerate, $V$ is
a symplectic vector space, so this type of quandle is called a
\textit{symplectic quandle}. If $\langle,\rangle$ is instead a 
symmetric bilinear form, then the subset 
$S=\{\mathbf{x}\in V \ :\ \langle \mathbf{x},\mathbf{x}\rangle\ne 0\}
\subset V$ is a quandle under 
\[\mathbf{x}\tr\mathbf{y}=
2\frac{\langle\mathbf{x},\mathbf{y}\rangle}{\langle\mathbf{y},
\mathbf{y}\rangle}\mathbf{y}-\mathbf{x}\]
called a \textit{Coxeter quandle}.

Another important example is the \textit{knot quandle} $Q(L)$ defined in 
\cite{J} which associates a quandle generator to every arc in a link diagram 
$L$ and a relation at every crossing. The elements of a knot quandle are 
equivalence classes of quandle words in the generators modulo the 
equivalence relation generated by the quandle axioms and the relations 
imposed by the crossings. 

\[\includegraphics{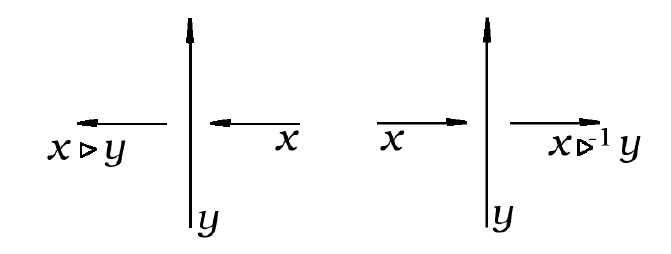}\]

\begin{example} \textup{The two-component link pictured has knot quandle
given by the listed presentation.}
\[\raisebox{-0.5in}{\includegraphics{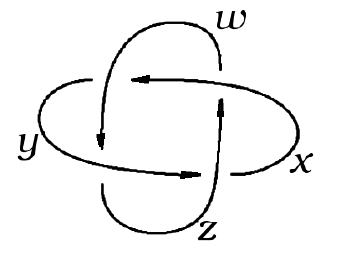}} 
\quad Q(L)=\langle x,y,z,w \ | \ x\tr z=y, \ z\tr y=w, y\tr w = x, 
w\tr x=z\rangle.\]
\end{example}

For symbolic computation, it is convenient to represent a finite
quandle $Q=\{x_1, x_2, \dots, x_n\}$ with an $n\times n$ matrix $M$ which
encodes the operation table of $Q$, i.e., $M_{ij}=k$ where
$x_k=x_i\tr x_j$. For example, the finite Alexander quandle 
$Q=\mathbb{Z}[t^{\pm 1}]/(2,t^2+1)=\{x_1=0,x_2=1,x_3=t,x_4=1+t\}$
has quandle operation matrix
\[M_Q=\left[\begin{array}{rrrr}
1 & 4 & 4 & 1 \\
3 & 2 & 2 & 3 \\
2 & 3 & 3 & 2 \\
4 & 1 & 1 & 4
\end{array}\right].\]

Quandles have been much studied in recent years; see \cite{J,FR,M,C3,Z,E}, etc.
Finite quandles are of particular interest since they can
be used to define an easily computable invariant of knots and links, the
\textit{quandle counting invariant} $|\mathrm{Hom}(Q(L),T)|$. Each 
quandle homomorphism $f\in\mathrm{Hom}(Q(L),T)$ can be pictured as a 
coloring of a link diagram representing the link $L$ with a quandle element 
$f(x)\in T$ attached to each arc $x$ satisfying the crossing condition
pictured above.

It is well known that quandle and biquandle counting invariants are stronger 
than fundamental group counting invariants (i.e., counting homomorphisms 
from a knot group to a finite group). For example, the \textit{Kishino virtual 
knot} below is distinguished from the unknot by biquandle counting invariants
despite having trivial knot group; see \cite{NV}. Similarly, the knot quandles
of the square knot $SK$ and granny knot $GK$ are known to be nonisomorphic
since their generalized knot groups $G_n$ derived from their knot quandles are 
distinguished by counting invariants for $n\ge 2$ in \cite{T}. We suspect 
that the 
corresponding conjugation quandle counting invariants should distinguish
the knots, though direct computational confirmation of this is difficult
due to the size of the groups involved.
However, the fact that $\pi_1(S^3\setminus GK)\cong \pi_1(S^3\setminus SK)$
implies that no group counting invariant $|\mathrm{Hom}(G,T)|$ where $G$ is 
the usual knot group can distinguish the square knot from the granny knot.
\[\begin{array}{ccc}
\includegraphics{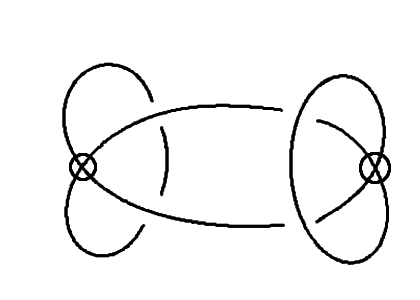} & \includegraphics{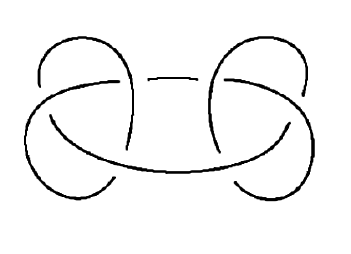} 
& \includegraphics{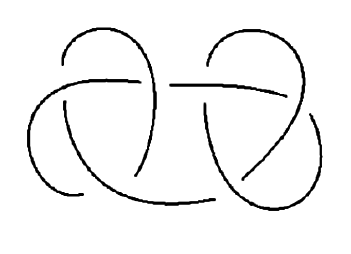} \\
\mathrm{Kishino} & \mathrm{Square} & \mathrm{Granny}
\end{array}
\]

In \cite{N} we have the following definition:

\begin{definition}\label{qp}
\textup{Let $Q$ be a finite quandle. For any element $x\in Q$, let}
\[ C(x) = \{ y\in Q \ : \ y\tr x = y\}
\quad \mathrm{and} 
\quad R(x) = \{ y\in Q \ : \ x\tr y = x \}\]
\textup{and set $r(x)=|R(x)|$ and $c(x)=|C(x)|$.
Then the \textit{quandle polynomial of $Q$}, $qp_Q(s,t)$, }
is \[qp_Q(s,t)=\sum_{x\in Q} s^{r(x)}t^{c(x)}.\] 
\end{definition}

An isomorphism $\phi:Q\to Q'$ induces bijections 
$\phi_r:R(x)\to R(\phi(x))$ and $\phi_c:C(x)\to C(\phi(x))$, so
$qp_Q(Q)=qp_Q(Q')$ and $qp_Q$ is an invariant of isomorphism
type for finite quandles.

We can now make our first new definition.

\begin{definition}\label{gqpdef}
\textup{Let $Q$ be a finite quandle. For any element $x\in Q$, let}
\[ C_n(x) = \{ y\in Q \ : \ y\tr^n x = y\}
\quad \mathrm{and \ let} 
\quad R_m(x) = \{ y\in Q \ : \ x\tr^m y = x \}\]
\textup{and set $r_m(x)=|R_m(x)|$ and $c_n(x)=|C_n(x)|$.
Then the \textit{$(m,n)$-quandle polynomial of $Q$}, $qp_{m,n}(Q)$, 
is} \[qp_{m,n}(Q)=\sum_{x\in Q} s^{r_m(x)}t^{c_n(x)}.\] 
\end{definition}

An isomorphism $\phi:Q\to Q'$ then induces for each $x\in Q$ bijections 
$\phi_{r,m}:R_m(x)\to R_m(\phi(x))$ and $\psi_{c,n}:C_n(x)\to C_n(\phi(x))$. 
It follows that $qp_{m,n}(Q)=qp_{m,n}(Q')$ and $qp_Q$ is an invariant of 
quandle isomorphism for finite quandles. 

\begin{example}\textup{From the definition, we have 
$qp_{0,0}(Q)=|Q|s^{|Q|}t^{|Q|}$ and $qp_{1,1}(Q)=qp_Q(Q)$. }
\end{example}

For any element $y\in Q$, the second quandle axiom implies that $y$
acts on $Q$ by a permutation $\rho_y\in S_{|Q|}$ (where $S_{|Q|}$ is the
symmetric group on $|Q|$ letters) given by the column
corresponding to $y$ in the matrix of $Q$. Let $\mathrm{ord}(\rho)$ denote the
order of $\rho$ in $S_{|Q|}$, i.e., the cardinality of the cyclic subgroup
of $S_{|Q|}$ generated by $\rho$. Then if 
$n\equiv n' \ \mathrm{mod}\ \mathrm{ord}(\rho_y)$, we have 
$x\tr^{n}y =x\tr^{n'} y$. It follows that for any finite quandle
$Q$, there are at most $N^2$ distinct generalized quandle polynomials
where $N=\mathrm{lcm}\{\mathrm{ord}(\rho_y)\ : \ y\in Q\}$. In particular,
to find all generalized quandle polynomials 
it suffices to consider the subset $\{qp_{m,n} \ | \ 0\le m,n \le N-1\}$
of the $\mathbb{Z}^2$-lattice of generalized quandle polynomials. For 
ease of comparison, we can write these
entries in an $N\times N$ matrix $M_{qp}(Q)$ whose $(i,j)$ entry is
$qp_{i-1,j-1}(Q)$, which we will call the \textit{generalized
quandle polynomial matrix} of $Q$. Both the size and the entries of this 
matrix are invariants of quandle isomorphism.

\begin{example}\textup{
The Alexander quandle 
$Q=\mathbb{Z}[t^{\pm 1}]/(2,t^2+1)=\{x_1=0,x_2=1,x_3=t,x_4=1+t\}$
has $\mathrm{ord}(\rho_y)=2$ for all $y\in Q$, so $N=2$; we compute the
generalized quandle polynomials 
$qp_{0,0}=4s^4t^4,\ qp_{0,1}=4s^2t^4, \ qp_{1,0}=4s^4t^2,$
and $qp_{1,1}=4s^2t^2$. Thus, $Q$ has generalized quandle polynomial
matrix}
\[M_{qp}(Q)=\left[\begin{array}{cc}
4s^4t^4 & 4s^2t^4 \\
4s^4t^2 & 4s^2t^2
\end{array}\right]. \]
\end{example}

\begin{example}\textup{A quandle is \textit{strongly connected} or 
\textit{Latin} if its operation matrix is a Latin square, that is, if
its rows as well as columns are permutations (see \cite{HN}). In \cite{N}, 
\texttt{Maple} computations showed that $qp_{1,1}(Q)$ distinguishes
all non-Latin quandles of cardinality up to 5, while Latin quandles 
always have $qp_{1,1}(Q)=|Q|st.$ Our \texttt{Maple} computations show that
of the three Latin quandles with five elements, two have the same
generalized quandle polynomial matrix while one has a different matrix. 
This shows that the generalized quandle polynomials contain additional 
information about quandle isomorphism type not contained in $qp_Q(s,t)$, 
while still not determining the quandle's isomorphism type for Latin quandles.
See Table 1.}
\end{example}

\begin{table}[!ht]
\begin{center}
\begin{tabular}{|c|c|} \hline
Quandle matrix & Generalized quandle polynomial matrix \\ \hline
& \\
$ \left[\begin{array}{ccccc} 
1 & 5 & 4 & 3 & 2 \\
3 & 2 & 1 & 5 & 4 \\
5 & 4 & 3 & 2 & 1 \\
2 & 1 & 5 & 4 & 3 \\
4 & 3 & 2 & 1 & 5
\end{array}\right]$ 
&
$ \left[\begin{array}{cccc} 
5s^5t^5 & 5st^5 & 5st^5 & 5st^5 \\
5s^5t   & 5st   & 5st   & 5st   \\
5s^5t   & 5st   & 5st   & 5st   \\
5s^5t   & 5st   & 5st   & 5st
\end{array}\right]$ \\
& \\
$\left[\begin{array}{ccccc} 
1 & 4 & 2 & 5 & 3 \\
4 & 2 & 5 & 3 & 1 \\
2 & 5 & 3 & 1 & 4 \\
5 & 3 & 1 & 4 & 2 \\
3 & 1 & 4 & 2 & 5
\end{array}\right]$ & 
$\left[\begin{array}{cccc} 
5s^5t^5 & 5st^5 & 5st^5 & 5st^5 \\
5s^5t   & 5st   & 5st   & 5st   \\
5s^5t   & 5st   & 5st   & 5st   \\
5s^5t   & 5st   & 5st   & 5st
\end{array}\right] $ \\
& \\
$\left[\begin{array}{ccccc} 
1 & 3 & 5 & 2 & 4 \\
5 & 2 & 4 & 1 & 3 \\
4 & 1 & 3 & 5 & 2 \\
3 & 5 & 2 & 4 & 1 \\
2 & 4 & 1 & 3 & 5 
\end{array}\right]$ & 
$\left[\begin{array}{cc} 
5s^5t^5 & 5st^5 \\
5s^5t   & 5st
\end{array}\right] 
$ \\ 
& \\ \hline
\end{tabular}
\end{center}
\caption{Generalized quandle polynomial matrices of Latin quandles of 
cardinality 5.}
\end{table}

\section{\large \textbf{Biquandle polynomials}}\label{bqp}

In the section we define the analog of generalized quandle
polynomials for finite biquandles. We start with the definition of a 
biquandle, also known as a type of \textit{switch} or 
\textit{Yang-Baxter Set}; see \cite{KR}.

\begin{definition}
\textup{A \textit{biquandle} is a set $B$ with four binary operations
$B\times B\to B$ denoted by
\[(a,b) \mapsto a^b, \ a^{\bar b}, \ a_b,\quad \mathrm{and} \quad a_{\bar b}\]
respectively, satisfying the following axioms:}
\begin{list}{}{}
\item[\textup{1.}]{\textup{For every pair of elements $a,b\in B$, we have}
\[
\mathrm{(i)} \ a=a^{b{\bar{b_a}}}, \quad
\mathrm{(ii)} \ b=b_{a{\bar{a^b}}}, \quad
\mathrm{(iii)} \ a=a^{\bar{b}b_{\bar a}}, \quad
\mathrm{and} \quad
\mathrm{(iv)} \ b=b_{\bar{a}a^{\bar b}}.
\]}
\item[\textup{2.}]{\textup{ Given elements $a,b\in B$, there are 
elements $x,y\in B$ such that}
\[
\mathrm{(i)} \ x=a^{b_{\bar x}}, \quad
\mathrm{ (ii) } \ a=x^{\bar b},\quad
\mathrm{ (iii)} \   b=b_{{\bar x}a}, \]\[
\mathrm{ (iv)} \  y=a^{\bar{b_y}},\quad
\mathrm{ (v) } \  a=y^b,   \quad \mathrm{and} \quad
\mathrm{ (vi)} \  b=b_{y{\bar a}}.
\]}
\item[\textup{3.}]{\textup{ For every triple $a,b,c \in B$ we have:}
\[
\mathrm{ (i)} \ a^{bc}=a^{c_bb^c}, \quad
\mathrm{(ii)} \ c_{ba} =c_{a^bb_a}, \quad
\mathrm{(iii)} \ (b_a)^{c_{a^b}}=(b^c)_{a^{c_b}}, \]
\[
\mathrm{(iv)} \ a^{{\bar b}{\bar c}}
=a^{\bar{c_{\bar b}}\bar{b^{\bar c}}},
\quad
\mathrm{(v)} \ c_{{\bar b}{\bar a}}
=c_{\bar{a^{\bar b}}{\bar {b_{\bar a}}}},
\quad \mathrm{and}  \quad
\mathrm{(vi)} \ (b_{\bar a})^{\bar{c_{{\bar {a^{\bar b}}}}}}
=(b^{\bar c})_{\bar{a^{\bar{c_{\bar b}}}}}.
\]}
\item[\textup{4.}]{\textup{Given an element $a\in B,$ there are 
\textit{unique} elements
$x,y\in B$ such that}
\[
\mathrm{ (i)} \   x=a_x, \quad
\mathrm{ (ii)} \  a=x^a, \quad
\mathrm{(iii)} \  y=a^{\bar y}, \quad \mathrm{and}  \quad
\mathrm{(iv) } \  a=y_{\bar a}.
\]}
\end{list}
\end{definition}

\[\includegraphics{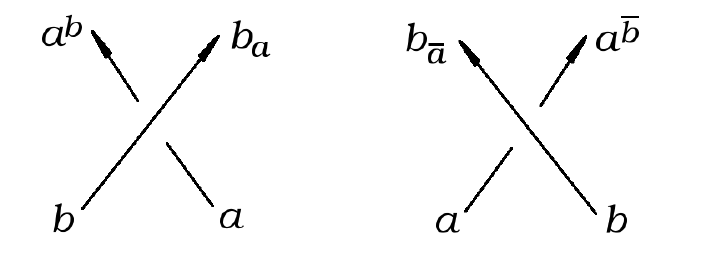}\]

The biquandle axioms come from dividing an oriented link diagram into
semi-arcs at every over and under crossing point; then both inbound 
semiarcs act on each other at both positive and negative crossings, for a
total of four binary operations $(a,b)\mapsto a^b, a^{\bar{b}}, a_b$
and $a_{\bar{b}}$. The axioms are then the result of transcribing a
minimal set of oriented Reidemeister moves. See \cite{KR} for more.

One standard example of a biquandle is any quandle, which is a biquandle under
\[a^b=a\tr b,\ \ a^{\bar b}=a\tr^{-1} b,\ \ a_b=a_{\bar b}=a\] as well as
under
\[a^b=a\tr^{-1} b,\\  a^{\bar b}=a\tr b,\ \ a_b=a_{\bar b}=a,\]
\[a_b=a\tr b,\ \ a_{\bar b}=a\tr^{-1} b,\ \ a^b=a^{\bar b}=a\] and
\[a_b=a\tr^{-1} b,\ \ a_{\bar b}=a\tr b,\ \ a^b=a^{\bar b}=a.\]

Anther standard example of a biquandle structure is any module over
$\mathbb{Z}[t^{\pm 1},s^{\pm 1}]$ with
\[a^b=ta+(1-st)b,\quad a^{\bar{b}}=t^{-1}a+(1-s^{-1}t^{-1})b,\quad
a_b=sa, \quad \mathrm{and} \quad a_{\bar{b}}=s^{-1}b;\]
biquandles of this type are called \textit{Alexander biquandles}.
For a concrete example, take $B=\mathbb{Z}_n$ and let $s,t\in B$ be any
two invertible elements. See \cite{KR} and \cite{LN}.

Much like with finite quandles, we can represent a finite biquandle
$B=\{x_1,\dots,x_n\}$ with
a block matrix encoding the operation tables of the four operations:
\[M_B=\left[\begin{array}{c|c} B^1 & B^2 \\ \hline 
B^3 & B^4 \end{array}\right] \quad 
B_{ij}^k=m \quad \mathrm{where} \quad x_m=\left\{\begin{array}{ll}
(x_i)^{\bar{(x_j)}} & k=1 \\
(x_i)^{(x_j)} & k=2 \\
(x_i)_{\bar{(x_j)}} & k=3 \\
(x_i)_{(x_j)} & k=4. \\
\end{array}\right.
\]

\begin{example}\label{abq1}
\textup{The Alexander biquandle $B=\mathbb{Z}_4$ with $s=3$ and $t=1$
has biquandle matrix}
\[ M_B=
\left[\begin{array}{rrrr|rrrr}
3 & 1 & 3 & 1 & 3 & 1 & 3 & 1 \\
4 & 2 & 4 & 2 & 4 & 2 & 4 & 2 \\
1 & 3 & 1 & 3 & 1 & 3 & 1 & 3 \\
2 & 4 & 2 & 4 & 2 & 4 & 2 & 4 \\ \hline
3 & 3 & 3 & 3 & 3 & 3 & 3 & 3 \\
2 & 2 & 2 & 2 & 2 & 2 & 2 & 2 \\
1 & 1 & 1 & 1 & 1 & 1 & 1 & 1 \\
4 & 4 & 4 & 4 & 4 & 4 & 4 & 4 \\
\end{array}\right] \quad \mathrm{where} 
\quad B=\{x_1=1,\ x_2=2,\ x_3=3,\ x_4=0\}.\]
\end{example}

In what follows, we will find it convenient to use the notation
\[\mathrm{op}_1(x,y)=x^{\bar{y}},\quad \mathrm{op}_2(x,y)=x^{y}, \quad
\mathrm{op}_3(x,y)=x_{\bar{y}}\quad \mathrm{and} \quad 
\mathrm{op}_4(x,y)=x_y,\] 
and as before $\mathrm{op}_i^n(x,y) = 
\mathrm{op}_i(\dots \mathrm{op}_i(\mathrm{op}_i(x,y),y)\dots y)$, where $n$ is
the number of ``$\mathrm{op_i}$''s.

We can extend the generalized quandle polynomials in the obvious way to 
obtain an invariant of biquandles up to isomorphism:

\begin{definition}
\textup{Let $B$ be a finite biquandle. For every $x\in B$, define}
\[C^i_n(x)=\{y\in B \ | \ \mathrm{op}_i^n(y,x)=y\}
\quad \mathrm{and} \quad 
R^i_m(x)=\{y\in B \ | \ \mathrm{op}_i^m(x,y)=x\}\]
\textup{where $m,n\in \mathbb{Z}$. Let $c^i_n(x)=|C^i_n(x)|$ and 
$r^i_m(x)=|R^i_m(x)|$ for 
$i=1,\dots,4$. Then the $(m,n)$ \textit{biquandle polynomial} of $B$ is}
\[bp_{m,n}(B)=\sum_{x\in B}
s_1^{r^1_m(x)}s_2^{r^2_m(x)}s_3^{r^3_m(x)}s_4^{r^4_m(x)}
t_1^{c^1_n(x)}t_2^{c^2_n(x)}t_3^{c^3_n(x)}t_4^{c^4_n(x)}.
\]
\end{definition}

\begin{example}
\textup{As in the quandle case,  we have
\[bp_{0,0}(B)=|B|s_1^{|B|}t_1^{|B|}s_2^{|B|}t_2^{|B|}
s_3^{|B|}t_3^{|B|}s_4^{|B|}t_4^{|B|}\]
for every finite biquandle. Also as in the quandle case, specializing
$s_i=t_i=1$ for all $i=1,\dots, 4$ yields $|B|$ for all $m,n\in \mathbb{Z}.$}
\end{example}

\begin{example}\label{exq}
\textup{Every quandle $Q$ is a biquandle with $a^b=a\tr b$, 
$a^{\bar{b}}=a\tr^{-1} b$ and $a_b=a_{\bar{b}}=a$. Specializing
$s_2=s, \ t_2=t$ and $s_i=t_i=1$ for $i=1,3,4$ in $bp_{m,n}(Q)$ yields
$qp_{m,n}(Q)$.}
\end{example}

\begin{example}\textup{
The Alexander biquandle $B=\mathbb{Z}_3$ with $s=2$, $t=1$ has biquandle 
matrix}
\[M_B=\left[\begin{array}{rrr|rrr}
3 & 2 & 1 & 3 & 2 & 1 \\
1 & 3 & 2 & 1 & 3 & 2 \\
2 & 1 & 3 & 2 & 1 & 3 \\ \hline
2 & 2 & 2 & 2 & 2 & 2 \\
1 & 1 & 1 & 1 & 1 & 1 \\
3 & 3 & 3 & 3 & 3 & 3
\end{array}\right].\]
\textup{We can compute $bp_{1,1}(B)$ by counting the number of times a
row number appears in each column and row in each of the operation block
matrices. Here we see that 
$bp_{1,1}(B)=2s_1s_2t_3t_4 + s_1t_1^3s_2t_2^3s_3^3t_3s_4^3t_4$.}
\end{example}

As before, a biquandle isomorphism $\phi:B\to B'$ induces bijections
$\phi_{i,n}:C^i_n(x)\to C^i_n(\phi(x))$ and 
$\psi_{i,m}:R^i_m(x)\to R^i_m(\phi(x))$ for each $i=1,\dots,4,$ and 
$n,m\in \mathbb{Z}$, so $B\cong B'$ implies $bp_{m,n}(B)=bp_{m,n}(B')$
and each $bp_{m,n}(B)$ is an invariant of biquandle isomorphism. 

The columns in a finite biquandle matrix, like those in a quandle matrix,
are permutations of $\{1, 2,\dots, |B|\}$. Let $\rho_1(y)\in S_{|B|},\
\rho_2(y)\in S_{|B|},\ \rho_3(y)\in S_{|B|},$ and $\rho_4(y)\in S_{|B|}$ be 
the permutations corresponding to the actions of $y$ on $B$ given by 
$\mathrm{op}_i(\underline{\ \ },y):B\to B$. Then as in the quandle case, the
fact that $n\equiv n' \ \mathrm{mod}\ \mathrm{ord}(\rho_i(y))$ implies 
$\mathrm{op}_i^n(x,y)=\mathrm{op}_i^{n'}(x,y)$ then implies that 
$bp_{m,n}(B)=bp_{m',n'}(B)$ if $n\equiv n'\ \mathrm{mod} \ N$ and 
$m\equiv m' \ \mathrm{mod}\ N$ where 
$N=\mathrm{lcm}\{\mathrm{ord}(\rho_i(y))\ :\ y\in B, \ i=1,\dots, 4\}.$
Hence, as before, there are at most $N^2$ distinct biquandle polynomials
for a biquandle $B$.

Thus, for every finite biquandle $B$, the \textit{biquandle
polynomial matrix} of $B$ is the $N\times N$ matrix whose $m,n$ entry is
$bp_{m,n}(B)$. Continuing with example \ref{exq}, if a biquandle $B$ is a 
quandle with $a\tr b=a^b, a\tr^{-1}b=a^{\bar{b}}$ and $a_b=a_{\bar{b}}=a$, 
then specializing $s_1=t_1=s_3=t_3=s_4=t_4=1$ and $s_2=s,\ t_2=t$ in the 
biquandle polynomial matrix of $B$ yields
the generalized quandle polynomial matrix $M_{qp}(B)$. 

\begin{example}\textup{The Alexander biquandle in example \ref{abq1}
has biquandle polynomial matrix}
\[ 
\left[\begin{array}{cc} 
4s_1^4t_1^4s_2^4t_2^4s_3^4t_3^4s_4^4t_4^4 &
2s_1^2t_1^4s_2^2t_2^4t_3^4t_4^4+2s_1^2t_1^4s_2^2t_2^4s_3^4t_3^4s_4^4t_4^4
\\
2s_1^4s_2^4s_3^4t_3^2s_4^4t_4^2+2s_1^4t_1^4s_2^4t_2^4s_3^4t_3^2s_4^4t_4^2 &
2s_1^2s_2^2t_3^2t_4^2+2s_1^2t_1^4s_2^2t_2^4s_3^4t_3^2s_4^4t_4^2
\end{array}\right] 
\]
\textup{according to our \texttt{Maple} computations.}
\end{example}

\texttt{Maple} code for computing quandle and biquandle polynomials
is available for download at \texttt{www.esotericka.org/quandles}. 
Computations with this code reveal that all isomorphism classes of
biquandles with up to four elements are distinguished by $bp_{1,1}(B)$ alone,
using the list of biquandle isomorphism classes from \cite{NV}.\footnote{Note
that the published list contains a few small typographical errors, but
one can regenerate the correct list with the \texttt{Maple} code.}

\section{\large \textbf{Link invariants from generalized quandle polynomials}}
\label{inv}

In this section we extend the subquandle polynomial defined in \cite{N} 
to the $qp_{m,n}(Q)$ and $bp_{m,n}(B)$ settings and exhibit some examples 
of the resulting link invariants.

\begin{definition}\textup{
Let $Q$ be a finite quandle and $S\subset Q$ a subquandle. Then for any
$m,n\in \mathbb{Z}$, the \textit{generalized subquandle polynomial} is}
\[sqp_{m,n}(S\subset Q)=\sum_{x\in S}s^{r_m(x)}t^{c_n(x)}. \]
\textup{Similarly, for any finite biquandle $B$ with subbiquandle 
$S\subset B$ and $m,n\in \mathbb{Z}$ the \textit{subbiquandle polynomial} is}
\[sbp_{m,n}(S\subset B)=\sum_{x\in S}
s_1^{r^1_m(x)}s_2^{r^2_m(x)}s_3^{r^3_m(x)}s_4^{r^4_m(x)}
t_1^{c^1_n(x)}t_2^{c^2_n(x)}t_3^{c^3_n(x)}t_4^{c^4_n(x)}.
\]
\end{definition}

Thus, the subquandle and subbiquandle polynomials are the contributions
to the quandle and biquandle polynomials coming from the subquandle or
subbiquandle in question. The $(1,1)$ subquandle polynomial was shown in 
\cite{N} to encode information about how the subquandle $S$ is embedded in
$Q$; indeed, $sqp_{1,1}$ can distinguish isomorphic subquandles embedded
in different ways.

Since the image of a homomorphism of a knot quandle into a target quandle
$T$ is a subquandle of $T$, we can modify the quandle counting invariant
to obtain a multiset-valued  link invariant by counting the subquandle 
polynomial of $\mathrm{Im}(f)$ for each $f\in \mathrm{Hom}(Q(L),T)$. 
The cardinality of this multiset is then the usual counting invariant. 
This generalizes the specialized subquandle polynomial invariant which was 
shown in \cite{N} to distinguish some links which have the same quandle 
counting invariant.

\begin{definition}\label{def-sqp}
\textup{Let $L$ be a link, $T$ a finite quandle and $m,n\in \mathbb{Z}$. 
Then the multiset}
\[\Phi_{sqp_{m,n}}(L,T)=
\{sqp_{m,n}(\mathrm{Im}(f)\subset T)\ :\ f\in\mathrm{Hom}(Q(L),T)\}\]
\textup{is the \textit{$(m,n)$-subquandle polynomial invariant} of
$L$ with respect to $T$. We can rewrite the multiset in a polynomial-style 
form by converting the multiset elements to exponents of a dummy variable 
$q$ and converting their multiplicities to coefficients:}
\[\phi_{sqp_{m,n}}(L,T)
=\sum_{f\in \mathrm{Hom}(Q(L),T)} q^{sqp_{m,n}(\mathrm{Im}(f)\subset T)}.\]

\textup{If $T$ is a finite biquandle, we similarly define the
\textit{$(m,n)$-subbiquandle polynomial invariant} to be the multiset}
\[\Phi_{sbp_{m,n}}(L,T)
=\{sbp_{m,n}(\mathrm{Im}(f)\subset T)\ :\ f\in\mathrm{Hom}(B(L),T)\}\]
\textup{or in polynomial form}
\[\phi_{sbp_{m,n}}(L,T)
=\sum_{f\in \mathrm{Hom}(B(L),T)} q^{sbp_{m,n}(\mathrm{Im}(f)\subset T)}.\]
\end{definition}

Collecting all the of the subquandle or subbiquandle polynomial
invariants into an $N\times N$ matrix whose $(m,n)$ entry is 
$sqp_{m,n}(\mathrm{Im}(f)\subset T)$ or $sbp_{m,n}(\mathrm{Im}(f)\subset T)$
yields an invariant of links which includes information from all of the 
subquandle or subbiquandle polynomials. Specializing 
$s_i=t_i=0$ for $i=1,\dots,4$ or specializing $q=1$ in any entry of the 
matrix yields the appropriate counting invariant.

\begin{example}

\[\raisebox{-0.5in}{\includegraphics{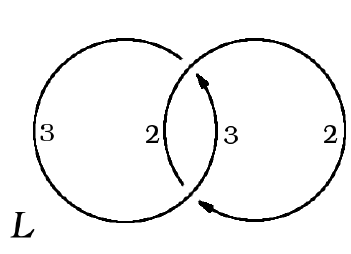}} \quad
M_T=\left[\begin{array}{ccc|ccc}
1 & 1 & 1 & 1 & 1 & 1 \\
3 & 2 & 2 & 3 & 2 & 2 \\
2 & 3 & 3 & 2 & 3 & 3 \\ \hline
1 & 1 & 1 & 1 & 1 & 1 \\
3 & 2 & 2 & 3 & 2 & 2 \\
2 & 3 & 3 & 2 & 3 & 3 \\
\end{array}
\right]\]

\textup{The Hopf link $L$ has nine biquandle colorings by the three-element 
biquandle $T$ below. In particular, the pictured coloring has image
subbiquandle $\mathrm{Im}(f)=\{2,3\}\subset T.$ The element $2\in T$
has $r^i_1=2$ and $c^i_1=3$ for $i=1,2,3,4$, so $2\in \mathrm{Im}(f)$
contributes $s_1^2s_2^2s_3^2s_4^2t_1^3t_2^3t_3^3t_4^3$ to the exponent 
of $q$ for this homomorphism. Indeed, the contribution from $3\in T$ is the 
same,
so we have a contribution of $q^{2s_1^2s_2^2s_3^2s_4^2t_1^3t_2^3t_3^3t_4^3}$
from this homomorphism. Repeating with the other homomorphisms, we have}
\begin{eqnarray*} \phi_{sbp_{1,1}}(L,T) & = &
q^{s_1^3t_1s_2^3t_2s_3^3t_3s_4^3t_4}
+4q^{s_1^3t_1s_2^3t_2s_3^3t_3s_4^3t_4
+2s_1^2t_1^3s_2^2t_2^3s_3^2t_3^3s_4^2t_4^3} \\
& & +2q^{s_1^2t_1^3s_2^2t_2^3s_3^2t_3^3s_4^2t_4^3}
+2q^{2s_1^2t_1^3s_2^2t_2^3s_3^2t_3^3s_4^2t_4^3}.\end{eqnarray*}
\textup{Since the least common multiple of the orders of the columns
of $T$ is 2, the full invariant is a $2\times 2$ matrix, of which the
above value is one entry.}
\end{example}

If $K$ is a single-component link, that is, a knot, then the knot quandle of
$K$ is connected, and hence the image of
any quandle homomorphism $f:Q(K)\to T$ must lie inside a single
orbit subquandle of $T$. In particular, we have 
\[|\mathrm{Hom}(Q(K),T)|=\left|\mathrm{Hom}\left(Q(K),
\bigcup_{i=1}^n T_i\right)\right|
=\sum_{i=1}^n|\mathrm{Hom}(Q(K),T_i)|\]
where $T_i$ are the orbit subquandles of $T$. In practice, this has meant
that multi-orbit quandles have been largely ignored in favor of single-orbit
(``connected'') quandles.
The next example demonstrates that by using generalized subquandle 
polynomials, non-connected quandles can still be used to distinguish knots 
whose counting invariants are the same.

\begin{example}
\textup{The quandle $T'$ with operation matrix}
\[M_{T'}=
\left[\begin{array}{cccccccccc} 
1 & 3 & 5 & 2 & 4 & 3 & 1 & 4 & 2 & 5 \\
5 & 2 & 4 & 1 & 3 & 5 & 3 & 1 & 4 & 2 \\
4 & 1 & 3 & 5 & 2 & 2 & 5 & 3 & 1 & 4 \\
3 & 5 & 2 & 4 & 1 & 4 & 2 & 5 & 3 & 1 \\
2 & 4 & 1 & 3 & 5 & 1 & 4 & 2 & 5 & 3 \\
8 & 9 & 10 & 6 & 7 & 6 & 10 & 9 & 8 & 7 \\
7 & 8 & 9 & 10 & 6 & 8 & 7 & 6 & 10 & 9 \\
6 & 7 & 8 & 9 & 10 & 10 & 9 & 8 & 7 & 6 \\
10 & 6 & 7 & 8 & 9 & 7 & 6 & 10 & 9 & 8 \\
9 & 10 & 6 & 7 & 8 & 9 & 8 & 7 & 6 & 10 
\end{array}\right]\]
\textup{has two 5-element orbit subquandles. The two knots $5_1$ and $6_1$
pictured below both have counting invariant $|\mathrm{Hom}(Q(5_1),T')|=30=
|\mathrm{Hom}(Q(6_1),T')|$ with respect to $T'$. Indeed, the two knots 
have the same $\phi_{sqp_{1,1}}(K,T')$ value. However, the generalized
subquandle polynomial invariants with $m=2$ distinguish the knots, detecting 
the fact that the sets of homomorphisms are different despite having the
same cardinality.}
\[
\includegraphics{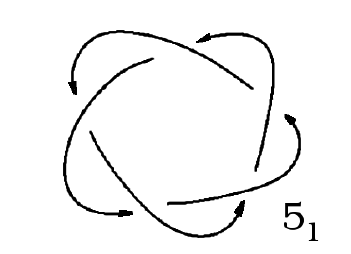} \quad \quad 
\includegraphics{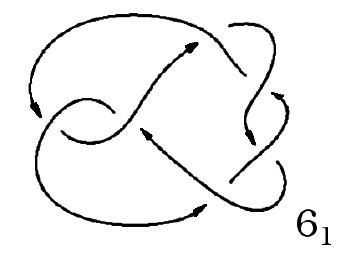}\]
\begin{center}
\begin{tabular}{|c|cc|} \hline
$n$ & $\phi_{sqp_{2,n}}(5_1,T')$ & $\phi_{sqp_{2,n}}(6_1,T')$ \\ \hline
& & \\
$0$ & $5q^{s^{10}t^{10}}+20q^{5s^{10}t^{10}}+5q^{s^{10}t^2}$ &
$5q^{s^{10}t^2}+20q^{5s^{10}t^2}+5q^{s^{10}t^{10}}$ \\
$1$ & $5q^{s^2t^{10}}+20q^{5s^2t^{10}}+5q^{s^2t^2}$ &
$5q^{s^2t^2}+20q^{5s^2t^2}+5q^{s^2t^{10}}$ \\
$2$ & $5q^{s^6t^{10}}+20q^{5s^6t^{10}}+5q^{s^6t^2}$ &
$5q^{s^6t^2}+20q^{5s^6t^2}+5q^{s^6t^{10}}$ \\
$3$ & $5q^{s^2t^{10}}+20q^{5s^2t^{10}}+5q^{s^2t^2}$ &
$5q^{s^2t^2}+20q^{5s^2t^2}+5q^{s^2t^{10}}$ \\
& & \\ \hline
\end{tabular}
\end{center}

\end{example}

\section{\large \textbf{Questions for future research}} \label{questions}

The quandle and biquandle polynomials as currently defined only make
sense for finite quandles. Consequentially, to use these polynomials 
for defining invariants of knot and link quandles and biquandles, which
are typically infinite, we must first convert to finite quandles in some way.
If a version of the quandle polynomial could be be defined for arbitrary 
quandles, or perhaps just for finitely generated quandles such as knot 
quandles, we might use such a polynomial (or series?) to obtain link 
invariants more directly. 

The only example so far of two non-isomorphic quandles or biquandles with
the same (bi)quandle polynomial matrix is the pair of Latin quandles
of order 5 from table 1. Are there examples of non-isomorphic non-Latin 
quandles or non-quandle biquandles which are not distinguished by their 
polynomial matrices?

If two knots or links have distinct quandles or biquandles, let the 
\textit{sub(bi)quandle polynomial matrix index} of the pair be the 
cardinality of the smallest finite (bi)quandle whose polynomial matrix 
invariant distinguishes the pair, or $\infty$ if there is no such
finite (bi)quandle. Is there a pair of knots or links whose sub(bi)quandle 
polynomial matrix index is infinite?

Each of the entries in a quandle or biquandle polynomial matrix has total
coefficient equal to the cardinality of the (bi)quandle. What other
relationships, if any, exist among the entries in a (bi)quandle polynomial
matrix? In particular, what is the minimal subset of the entries which
determines the other entries? 

\section{\large \textbf{Acknowledgements}}

The author would like to thank Jozef Przytycki for his comments and 
observations at the Knots in Washington conference which directly inspired
definition \ref{gqpdef}, as well as Thao-Nhi Luu, whose conversations with 
the author also influenced the direction of this paper, and the referee, whose
comments and observations improved the paper.

\begin{tabular}{l}
\textit{Department of Mathematics, Pomona College} \\
\textit{610 N. College Ave} \\
\textit{Claremont, CA 91711} \\
\textit{e-mail: knots@esotericka.org}
\end{tabular}

\end{document}